\documentclass[12pt]{article}
\usepackage{mathrsfs}
\usepackage{amssymb}
\usepackage{amsmath}
\usepackage{hyperref}
\usepackage{marvosym}
\usepackage{pict2e}
\usepackage{geometry}
\usepackage{color}
\def\q{\hfill\rule{1ex}{1ex}}
\def\0{\emptyset}

\def\p{{\bf Proof.}\quad}
\def\q{\hfill\rule{1ex}{1ex}}

\def\n{\noindent}

\newtheorem{theo}{Theorem}[section]
\newtheorem{coro}[theo]{Corollary}
\newtheorem{conj}[theo]{Conjecture}

\begin{document}

\title{\bf Chords of longest cycles passing through a specified small set}
\author{Haidong Wu~~~~~Shunzhe Zhang\\}

\author{Haidong Wu\footnote{E-mail: hwu@olemiss.edu}
~~~~~Shunzhe Zhang\footnote{E-mail: szhang12@go.olemiss.edu}\\
Department of Mathematics\\
The University of Mississippi\\
University, MS 38677, U.S.A\\
}

\date{}
\maketitle
\baselineskip 16.7pt \setcounter{page}{1}

\begin{abstract}

A long-standing conjecture of Thomassen says that every longest cycle of a $3$-connected graph has a chord. Thomassen (2018) proved that if $G$ is $2$-connected and cubic, then any longest cycle must have a chord. He also showed that if $G$ is a $3$-connected graph with minimum degree at least $4$, then some of the longest cycles in $G$ must have a chord. Zhang (1987) proved that if $G$ is a $3$-connected simple planar graph which is 3-regular or has minimum degree at least $4$, then every longest cycle of $G$ must have a chord. Recently, Li and Liu showed that if $G$ is a $2$-connected cubic graph and $x, y$ are two distinct vertices of $G$, then every longest $(x,y)$-path of $G$ contains at least one internal vertex whose neighbors are all in the path. In this paper, we study chords of longest cycles passing through a specified small set and generalize Thomassen's and Zhang's above results by proving the following results. (i) Let $G$ be a $2$-connected cubic graph and $S$ be a specified set consisting of an edge plus a vertex. Then every longest cycle of $G$ containing $S$ must have a chord. (ii) Let $G$ be a $3$-connected graph with minimum degree at least $4$ and $e$ be a specified edge of $G$. Then some longest cycle of $G$ containing $e$ must have a chord. (iii) Let $G$ be a $3$-connected planar graph with minimum degree at least $4$. Suppose $S$ is a specified set consisting of either three vertices or an edge plus a vertex. Then every longest cycle of $G$ containing $S$ must have a chord. We also extend the above-mentioned result of Li and Liu for $2$-connected cubic graphs. 

\vskip 0.1cm


\end{abstract}

\section{Introduction}

In 1976, Thomassen raised the following famous conjecture, see \cite{Bakr2010}, \cite{Bo1995}, \cite{Bo2014}.

\begin{conj}\label{conj1.1} {\rm (Thomassen~\cite{Al1985},\cite{Th1989})}
\label{originalconj}
Every longest cycle of a $3$-connected graph has a chord.
\end{conj}

Although the general conjecture remains unsolved, many partial results have been discovered. The first result regarding planar graphs was due to Zhang~\cite{Zh1987} in 1987, who showed the conjecture holds for 3-connected planar graphs which are either cubic or with a minimum degree at least four. 

\begin{theo}\label{theo1.2} {\rm (Zhang~\cite{Zh1987})}
Let $G$ be a $3$-connected planar graph which is either cubic or with minimum degree at least $4$. Then any longest cycle of $G$ must have a chord. \q
\end{theo}

Thomassen~\cite{Th1997} showed that Conjecture \ref{originalconj} holds when $G$ is cubic. In 2018, Thomassen ~\cite{Th2018} further extended his result to $2$-connected cubic graphs. He also showed that if $G$ is $3$-connected and has a minimum degree of at least $4$, then some longest cycle in $G$ has a chord. 

\begin{theo}\label{theo1.3} {\rm (Thomassen~\cite{Th2018})}
Every longest cycle in a $2$-connected cubic graph has a chord.
\end{theo}

\begin{theo}\label{theo1.4} {\rm (Thomassen~\cite{Th2018})}
Let $G$ be a $3$-connected graph with minimum degree at least $4$. Then some longest cycle of $G$ has a chord.
\end{theo}

Recently, Li and Liu generalized Theorem \ref{theo1.3} as follows.

\begin{theo}\label{theo1.5} {\rm (Li and Liu~\cite{Li2024})}
if $G$ is a $2$-connected cubic graph and $x, y$ are two distinct vertices of $G$, then every longest $(x,y)$-path of $G$ contains at least one internal vertex whose three neighbors are all in the path.
\end{theo}

Considering sparse graphs more generally, Li and Zhang~\cite{LiZh20031} verified the conjecture for graphs embedded in the projective plane with minimum degree at least $4$, and Li and Zhang~\cite{LiZh20032} verified the conjecture when G is $4$-connected and embedded in a torus or a Klein bottle. Kawarabayashi et al. ~\cite{Ka2007} verified the conjecture for locally 4-connected planar graphs, and Birmele ~\cite{Bi2008} verified the conjecture for $3$-connected graphs with no $K_{3,3}$-minor. Harvey~\cite{Ha2017} proved that  every longest cycle has a chord for graphs with a large minimum degree.

Conjecture \ref{originalconj} concerns the global structure in 3-connected graphs. In this paper, we consider the local structure in 3-connected graphs by studying longest cycles containing a specified small set. In 2005, Gu, Jia and Wu \cite{Wu2005} proposed the following conjecture.

\begin{conj}\label{conj1.4}{\rm (Gu, Jia and Wu~\cite{Wu2005})}
Let $G$ be a $3$-connected graph and $e$ be an edge of $G$. Then among all the cycles containing $e$, each longest such cycle has a chord.
\end{conj}

Inspired by Conjecture 1.6, we study longest cycles containing a specified small set and generalize Theorems 1.2-1.5 as follows. Our results also partially verify Conjecture 1.6. Our first result generalizes Theorem \ref{theo1.5}.

\begin{theo}\label{main1} 
Let $G$ be a $2$-connected cubic graph and $x,y,z$ be three distinct vertices in $G$. Then every longest $(x,y)$-path of $G$ containing $z$ has at least one internal vertex that has three neighbors in the path.
\end{theo}

The following corollary generalizes Theorem \ref{theo1.3}. 

\begin{coro}\label{main1.1}
Let $G$ be a $2$-connected cubic graph and $S$ be a specified set consisting of an edge plus a vertex. Then every longest cycle of $G$ containing $S$ must have a chord.
\end{coro}

\p Let $S=\{e, z\}$ be a specified set containing edge $e=xy$ and vertex $z$, and $C$ be a longest cycle of $G$ containing $S$. Then $P:=C-xy$ must be a longest $(x,y)$-path containing $z$. By Theorem \ref{main1}, $P$ contains an internal vertex that has three neighbors in $P$. Thus $C$ must have a chord. \q

\begin{coro}\label{main1.2}
Let $G$ be a $2$-connected cubic graph and $S$ be a specified set consisting of at most two vertices. Then every longest cycle of $G$ containing $S$ must have a chord.
\end{coro}

Our next result generalizes Theorem \ref{theo1.4}. 

\begin{theo}\label{main2} 
\label{main2}
Let $G$ be a $3$-connected graph with minimum degree at least $4$ and $e$ be a specified edge of $G$. Then some longest cycle of $G$ containing $e$ must have a chord.
\end{theo}

The following result generalizes the corresponding theorem of Zhang~\cite{Zh1987}.

\begin{theo}\label{main3} 
\label{main3}
Let $G$ be a $3$-connected planar graph with minimum degree at least $4$. Suppose $S$ is a specified set consisting of either three vertices or an edge plus a vertex. Then every longest cycle of $G$ containing $S$ must have a chord.
\end{theo}

By Corollary \ref{main1.1} and Theorem \ref{main3}, we obtain the following generalization of Theorem \ref{theo1.2}.

\begin{coro}\label{main3.1} 
Let $G$ be a $3$-connected planar graph which is either cubic or with minimum degree at least $4$ and $S$ be a set consisting of an edge plus a vertex or a set with at most two vertices. Then every longest cycle of $G$ containing $S$ must have a chord.
\end{coro}

Proofs of our main results will be given in Sections 3, 4 and 5. While using the techniques developed by Thomassen \cite{Th2018}, we also use some new ideas to prove our first two main results.  Throughout this paper, we consider finite simple graphs with no loops or multiple edges. We use Bondy and Murty~\cite{Bo2008} for terminology and notation not defined here. For a graph $G = (V, E)$, we use $V(G)$ and $E(G)$ to denote the vertex set and the edge set, respectively. A \emph{chord} of a cycle $C$ is an edge between two non-consecutive vertices of the cycle. We use $ch(C)$ to denote the number of chords of a cycle $C$. For a vertex $v\in V(G)$, $E(v)$ is the set of all edges incident with $v$ and $N_G(v)$ is the set of all neighbors of $v$ and $d_G(v) = |N_G(v)|$. If $S\subseteq V(G)$, then $G[S]$ denotes the subgraph induced by $S$. For a subgraph $H$ of $G$, let $G-H:= G[V(G)\setminus V(H)]$. An \emph{independent set} in a graph is a set of vertices no two of which are adjacent. A cycle in a graph $G$ is called \emph{Hamiltonian} if it contains each vertex of $G$. Let $C$ be a cycle/path with an arbitrary orientation. For $u\in V(C)$, we use $u^+$ to denote the successor of $u$ and $u^-$ its predecessor. If $u, v\in V(C)$, we use $uCv$ and {$v\overline{C}u$ to denote the subpath $uu^+\ldots v^-v$ of $C$ and the same subpath in reverse order, respectively. 

A face is said to be \emph{incident} with the vertices and edges in its boundary, and two faces are \emph{adjacent} if their boundaries have an edge in common. We denote the \emph{boundary} of a face $F$ by $\partial(F)$. Let $P$ be a subgraph of a graph $G$. A \emph{$P$-bridge} of $G$ is either an edge of $G-E(P)$ with both ends on $P$ or a subgraph of $G$ induced by the edges in a component of $G-V(P)$ and all edges from that component to $P$. For a $P$-bridge $B$ of $G$, the vertices in $V(B)\cap V(P)$ are the \emph{attachments} of $B$ on $P$. Define $A(B)=V(B)\cap V(P)$ and $I(B) = V(B)-A(B)$. Two open sets into which a simple closed curve $C$ partitions the plane are called the interior and exterior of $C$, respectively, which will be denoted by $int(C)$ and $ext(C)$. As in ~\cite{Zh1987}, a vertex set $T\subseteq V(G)$ is called a \emph{separating $k$-vertex-cut} with respect to $C$ if there is a partition set $V_1$ and $V_2$ of $G-T$ with $|T|=k$ such that $V(C)\cap V_1\neq \emptyset$ and $V(C)\cap V_2\neq \emptyset$. A cycle $C$ is called \emph{separable} if there is a separating $3$-vertex-cut with respect to $C$.

\section{Lemmas}

In this section, we prove a preliminary result which will be used in our proofs. This result is similar to Theorem 2 of Thomassen $\cite{Th2018}$. 

\n{\bf Lemma 2.1} {\em ~Let $G$ be a connected graph and $C=v_1v_2\ldots v_nv_1$ be a cycle of $G$ containing an odd degree vertex $v_i$ for some $i\in\{3,n\}$ such that all vertices in $G-V(C)$ have even degree in $G$. Suppose that no two vertices of even degree are adjacent except that the edge $v_1v_2$ is allowed. Then $G$ has a cycle $C'$ distinct from $C$ such that $C'$ contains all vertices of odd degree and the edge $v_1v_2$. 
}

\p Let $C=v_1v_2\ldots v_nv_1$ be a cycle of $G$ containing an odd degree vertex $v_i$ for some $i\in\{3,n\}$ such that all vertices in $G-V(C)$ have even degree in $G$. Choose a cycle $C$ such that no vertex in $G-V(C)$ is adjacent to two consecutive vertices of $C$ except $v_1$ and $v_2$. Otherwise, we can choose a longer cycle that contains all vertices of $C$ and $v_1v_2$ and one more vertex instead, which implies the conclusion is true. Note that $v_i$  is an odd degree vertex for some $i\in\{3,n\}$. Then the path $P^*=C-v_2v_3=v_2v_1\ldots v_3$ or $P^*=C-v_1v_n=v_1v_2\ldots v_n$ starts with the edge $v_1v_2$, contains all vertices of odd degree, and ends with a vertex $v_3$ or $v_n$ of odd degree. Without loss of generality, assume that $v_n$ has odd degree in $G$ and $P^*=C-\{v_1v_n\}=v_1v_2\ldots v_n$.

Applying the lollipop method, we consider an auxiliary graph $H$. A vertex in $H$ is a path $P$ in $G$ which starts with the edge $v_1v_2$, contains all vertices of odd degree, and ends with a vertex of odd degree. Thus $P^*\in V(H)$. Let $P_1=x_1x_2\ldots x_n$ and $P_2=y_1y_2\ldots y_n$ be two vertices in $H$, where $x_1=y_1=v_1$, $x_2=y_2=v_2$ and $x_n,y_n$ are two vertices of odd degree. Note that no two vertices of even degree are adjacent except that the edge $v_1v_2$ is allowed. $P_1$ and $P_2$ are adjacent in $H$ if $x_sx_n\in E(G)$ or $x_sux_n$ is a path with $u\in V(G)-V(P_1)$ for some $2\leq s\leq n-2$ and
 \begin{equation*}
P_2=\begin{cases}
~x_1x_2P_1x_sx_n\overline{P_1}x_{s+1},         &\mbox{if}~d_G(x_{s+1}) \mbox{~is~odd} ~\mbox{and}~x_sx_n\in E(G);       \\
~x_1x_2P_1x_sux_n\overline{P_1}x_{s+1},        &\mbox{if}~d_G(x_{s+1}) \mbox{~is~odd} ~\mbox{and}~x_sux_n \mbox{is a path};       \\
~x_1x_2P_1x_sx_n\overline{P_1}x_{s+2},         &\mbox{if}~d_G(x_{s+1}) \mbox{~is~even} ~\mbox{and}~x_sx_n\in E(G).      \\
~x_1x_2P_1x_sux_n\overline{P_1}x_{s+2},        &\mbox{if}~d_G(x_{s+1}) \mbox{~is~even}~\mbox{and}~x_sux_n \mbox{is a path}.  

\end{cases}
\end{equation*}
Let $U_P=\{u\in V(G)-V(P)\,|\,x_1ux_n~\mbox{is a path}\}$ where $P=x_1x_2\ldots x_n$ is a vertex in $H$. Since no two vertices of even degree are adjacent except that the edge $v_1v_2$ is allowed, we deduce that if $x_1x_n\notin E(G)$, then $d_H(P)=d_G(x_n)-(|U_P|+1)$ and if $x_1x_n\in E(G)$, then $d_H(P)=d_G(x_n)-(|U_P|+2)$. Note that $x_n$ has odd degree in $G$. It follows that $d_H(P)$ is odd if and only if either $x_1x_n\notin E(G)$ and $|U_P|$ is odd or $x_1x_n\in E(G)$ and $|U_P|$ is even. By the choice of $C$, $|U_{P^*}|=0$ and thus $P^*=C-\{v_1v_n\}=v_1v_2\ldots v_n$ has odd degree in $H$. Since the number of odd degree vertices in $H$ is even, there is another vertex $P'=v_1v_2u_3\ldots u_n$ in $H$ distinct from $P^*$ such that $d_H(P')$ is odd. Then either $v_1u_n\notin E(G)$ and $|U_{P'}|$ is odd or $v_1u_n\in E(G)$ and $|U_{P'}|$ is even. In either case, $G$ has a cycle $C'$ with either $C'=P'+\{v_1u,uu_n\}$, where $u\in V(G)-V(P')$ or $C'=P'+\{v_1u_n\}$ distinct from $C$ such that $C'$ contains all vertices of odd degree and $v_1v_2$.
\q

We will use the following results.

\n{\bf Theorem 2.2} {\rm (Thomassen~\cite{Th2018})} {\em ~Let $G$ be a cubic graph with a partition set $A$ and $B$ of $V(G)$ such that $G[A]$ is a matching $M$ and $G[B]$ is a matching $M'$ with $|A|=|B|=2k$. If $G$ has a cycle $C$ of length $3k$ which contains $M$ and precisely one end of each edge in $M'$. Then $G$ has a cycle of length larger than $3k$ containing $M$.}

\n{\bf Theorem 2.3} {\rm (Fleischner and Stiebitz~\cite{FlSt1992})}
{\em ~Let $n$ be a positive integer, and let $G$ be a $4$-regular graph on $3n$ vertices. Assume that $G$ has a decomposition into a Hamilton circuit and $n$ pairwise vertex disjoint triangles. Then $\chi(G)=3$.
}

\n{\bf Tutte's Lemma} {\rm (in Ore~\cite{Or1967})} {\em ~Let $G$ be a $2$-connected planar graph, let $e$ be an edge of $G$, let $F$ and $F'$ be the two faces incident with $e$ and let $e'$ be an edge on the boundary of $F$ and adjacent with $e$. Then there is a cycle $C$ of $G $ such that for any bridge $B$ of $C$, 

$(i)$ $e,e'\in E(C)$;

$(ii)$ $|A(B)|\leq 3$;

$(iii)$ If $E(B)\cap (\partial(F)\cup \partial(F'))\neq \emptyset$, then $|A(B)|\leq 2$. }

\n{\bf The Jordan Curve Theorem} {\rm (see Bondy and Murty~\cite{Bo2008})} {\em ~Any simple closed curve $C$ in the plane partitions the rest of the plane into two disjoint arcwise-connected open sets.
}

\section{Chords of longest cycles containing a specified small set in $2$-connected cubic graph}

In this section, we prove our first main result, Theorem \ref{main1}.   For convenience, we restate the theorem here. 
\begin{theo}\label{re-main1}
Let $G$ be a $2$-connected cubic graph and $x,y,z$ be three distinct vertices in $G$. Then every longest $(x,y)$-path of $G$ containing $z$ has at least one internal vertex that has three neighbors in the path.
\end{theo}

\p Let $G$ be a $2$-connected cubic graph and $x,y,z$ be three distinct vertices in $G$. Then there exists a $(x,y)$-path of $G$ containing $z$. In fact, if $P_1$ is a $(x,y)$-path of $G$ and $z\notin V(P_1)$, then there are two internally disjoint paths $zP_2a$ and $zP_3b$ from $z$ to $P_1$ as $G$ is 2-connected. Assume that $x,a,b,y$ on $P_1$ appear in the order $x,a,b,y$. Then $xP_1aP_2zP_3bP_1y$ is a $(x,y)$-path  of $G$ containing $z$. Suppose that $P$ is a longest $(x,y)$-path of $G$ containing $z$ but every internal vertex of $P$ has only two neighbors in $P$. Let $H_1,\ldots, H_l$ be all components of $G-V(P)$. Since $G$ is cubic, every internal vertex in $P$ has exactly one neighbor in $H_i$ and every vertex in $\{x,y\}$ has exactly one neighbor in $H_i$ if $xy\in E(G)$. Otherwise, every vertex in $\{x,y\}$ has exactly two neighbors in $G-V(P)$. Let $P=xu\cdots y$ with $z\in V(P)\setminus \{x,y\}$.

We construct a graph $G_1$ from $G$ as follows:

(1) If $H$ is a component of $G-V(P)$ joined to at least three vertices of $P$, then we contract $H$ to a single vertex. If this vertex is not adjacent to $x$, then we call this vertex as a \emph{pleasant vertex}. Otherwise, we call it as an \emph{unpleasant vertex}.

(2) If $H$ is a component of $G-V(P)$ joined to only two vertices $a,b$ of $P$, then we replace $H$ by an edge $ab$. If this edge is not incident with $x$, then we call this edge as a \emph{pleasant edge}. Otherwise, we call it as an \emph{unpleasant edge}.

Then $P$ is also a longest $(x,y)$-path of $G_1$ containing $z$. Since $G$ is a $2$-connected graph, $G_1$ is also a $2$-connected graph. Let $I=\{u_1,u_2,\ldots, u_p\}$ and $J=\{v_1,v_2,\ldots, v_m\}$ be the set of all pleasant vertices and unpleasant vertices, respectively. Let $D=\{e_1,e_2,\ldots, e_q\}$ and $F=\{f_1,f_2,\ldots, f_n\}$ be the set of all pleasant edges and unpleasant edges. Then $|I|+|J|+|D|+|F|=l$. Since $G$ is cubic, we have $d_{G_1}(x)=3$ and thus $0\leq m,n\leq 2$. Let $u\in N_P(H_i)$ for some $i\in \{1,\ldots, l\}$. If $x\in N_P(H_i)$, then assume that $xx',uu'\in E(G)$ with $x',u'\in V(H_i)$ and $P'$ is a $(x',u')$-path in $H_i$. Then $P^*=xx'P'u'u\overline{P}y$ is a $(x,y)$-path containing $z$ longer than $P$, a contradiction. Thus $x\notin N_P(H_i)$. This implies that $I\cup D\neq \emptyset$. Let $N_P(u_i)=\{u_{i1},u_{i2},\ldots, u_{it_i}\}$ for $u_i\in I$ and $N_P(v_j)=\{v_{j1},v_{j2},\ldots, v_{js_j}\}$ for $v_j\in J$. Then $t_i,s_j\geq 3$.

Now we construct a graph $G_1^*$ from $G_1$ as follows:

($1^*$) We select three neighbors $u_{i1},u_{i2},u_{i3}$ of $u_i$ on $uPy$ for each $u_i\in I$ and add edges $u_{i1}u_{i2},u_{i1}u_{i3},u_{i2}u_{i3}$.

($2^*$) We delete all vertices in $I\cup J\cup D\cup F$ and the vertex $x$ and add an edge $uy$.

($3^*$) We contract all edges each of which has one end vertex of degree two on the cycle $uPy+\{uy\}$. 

Then $G_1^*$ is a 4-regular graph on $3p$ vertices and $G_1^*$ has a decomposition into a Hamiltonian cycle and $p$ pairwise vertex disjoint triangles. By Theorem 2.3, we have $\chi(G_1^*)=3$. Since each vertex in a triangle must have different colors, we can assign a vertex of each triangle into color red. Let $R$ be a red class such that $\{u,z\}\cap R=\emptyset$. It is possible that $y\in R$. Then $R$ is independent in $G_1^*$. Thus $R$ is also independent in $G_1$. Let $R=\{r_1,r_2,\ldots, r_p\}$. Then there is a bijection $\varphi_1: I\rightarrow R$ denoted by $\varphi_1(u_i)=r_i$ for $i\in \{1,2,\ldots,p\}$.

If $xy\notin E(G_1)$, then we add an edge $xy$ in $G_1$. With a slight abuse of notation, we call the resulting graph $G_1$.  Since $G$ is cubic, we have $d_{G_1}(b)=3$ for $b\in V(C)\setminus \{x,y\}$, $d_{G_1}(x)=d_{G_1}(y)=3$ if $xy\in E(G_1)$ and $d_{G_1}(x)=d_{G_1}(y)=4$ if $xy\notin E(G_1)$. 

Let $C=P+xy$. A cycle $C_1$ is called \emph{pleasant} in $G_1$ if it contains all vertices of $C$ and $xy$ except some vertices in $R$ as in $\cite{Th2018}$. Then $xy\in E(C)\cap E(C_1)$. Combining this with $\{u,z\}\cap R=\emptyset$, we have $\{u,x,y,z\}\subseteq V(C)\cap V(C_1)$. Without loss of generality, assume that $R_1=\{r_1,\ldots, r_k\}\subseteq R$ is the set of all vertices in $C$ but not in $C_1$. Then $u,y\notin R_1$. Let $I_1$ and $J_1$ be the set of all pleasant vertices and unpleasant vertices in $C_1$, respectively. Let $D_1$ and $F_1$ be the set of all pleasant edges and unpleasant edges in $C_1$, respectively. Let $a\in R_1$. Then $C_1$ must contain two vertices adjacent to $a$ on $C$. Let $b$ be a neighbor of $a$ on $C$. Then $C_1$ contains an edge $bc$ or a path $bdc$, where $d\in I_1\cup J_1$ and $bc\in D_1\cup F_1$. We say the edge $bc$ or the vertex $d$ \emph{dominates} $a$ as in \cite{Th2018}.

{\bf Claim} {\em If a cycle $C_1$ is pleasant in $G_1$, then $|E(C)|=|E(C_1)|$, $|I_1|=|R_1|$ and $|J_1|=|D_1|=|F_1|=0$. Each vertex in $R_1$ is dominated by one vertex or two vertices in $I_1$ and each vertex in $I_1$ dominates one vertex or two vertices in $R_1$. In particular, the vertex $a$ in $R_1$ is precisely dominated by one vertex in $I_1$ if and only if this vertex in $I_1$ precisely dominates the vertex $a$ in $R_1$.
}

\p Let $b$ and $b'$ be two neighbors of $a$ on $C$. Then the edge $bc$ or the vertex $d$ dominates $a$. Possibly, the edge $bc$ or the vertex $d$ dominates some neighbor of $c$ on $C$. Note that $xy\in E(C)\cap E(C_1)$ and $d_{G_1}(b)=3$ for $b\in V(C)\setminus \{x,y\}$. Then each element in $I_1\cup J_1\cup D_1\cup F_1$ dominates at most two vertices in $R_1$. Similarly, $C_1$ contains an edge $b'c'$ or a path $b'd'c'$, where $d'\in I_1\cup J_1$ and $b'c'\in D_1\cup F_1$. Then the edge $b'c'$ or the vertex $d'$ dominates $a$. Maybe $a$ is precisely dominated by one element when the edge $bb'$ or the vertex $d=d'$ dominates $a$. So each vertex $a$ in $R_1$ is dominated by one element or two elements in $I_1\cup J_1\cup D_1\cup F_1$. Recall that any element in $I_1\cup J_1\cup D_1\cup F_1$ dominates at most two vertices in $R_1$. If follows that $|I_1|+|J_1|+|D_1|+|F_1|\geq |R_1|$.

Note that $|E(C_1)|=|E(C)|-2|R_1|+2(|I_1|+|J_1|)+(|D_1|+|F_1|)$. We can extend the edge set of the cycle $C_1$ containing $xy$ and $z$ in $G_1$ to the edge set of a cycle $C^*$ containing $xy$ and $z$ in $G$ as follows: each pleasant and unpleasant edge in $C_1$ can be extend to a path with at least 3 edges in $G$ as $G$ is cubic and each pleasant and unpleasant vertex in $C_1$ can be extend to a path with at least 1 edge in $G$. Since $P$ is a longest $(x,y)$-path of $G$ containing $z$ and $|I_1|+|J_1|+|D_1|+|F_1|\geq |R_1|$, we have $C$ is a longest cycle containing $xy$ and $z$ in $G$ and thus
$|E(C)|\geq |E(C^*)|\geq |E(C_1)|=|E(C)|-2|R_1|+2(|I_1|+|J_1|)+3(|D_1|+|F_1|)\geq |E(C)|+(|D_1|+|F_1|)\geq |E(C)|$. This implies that $|E(C)|=|E(C^*)|=|E(C_1)|$, $|I_1|+|J_1|=|R_1|$ and $|D_1|=|F_1|=0$. Then $C_1$ contains no pleasant and unpleasant edges and each vertex in $R_1$ is dominated by one vertex or two vertices in $I_1\cup J_1$ and each vertex in $I_1\cup J_1$ dominates one vertex or two vertices in $R_1$. In particular, the vertex $a$ in $R_1$ is precisely dominated by one vertex in $I_1\cup J_1$ if and only if this vertex in $I_1\cup J_1$ precisely dominates the vertex $a$ in $R_1$. Recall that $u,y\notin R_1$ and $ux,xy\in E(C)$. Then $|J_1|=0$.
\q

Now we construct a graph $G_2$ from $G_1$ as follows: 

($1'$) For any vertex $u_i\in I$ and its corresponding red vertex $r_i\in R$ with $\varphi_1(u_i)=r_i$, we contract the edge $u_ir_i$ to a vertex $r^*_i$.

($2'$) For any vertex $v_j\in J$, we contract the edge $v_jx$ to a vertex that is also denoted by $x$.

Let $R^*=\{r^*_1,r^*_2,\ldots, r^*_p\}$. Then there is a bijection $\varphi_2: R\rightarrow R^*$ denoted by $\varphi_2(r_i)=r^*_i$. So $G_2$ has a Hamiltonian cycle $C=xu\cdots yx$ and $d_{G_2}(b)=3$ for any $b\in V(C)\setminus R^*\cup \{x,y\}$. Since $R$ is independent in $G_1$, we have $R^*$ is also independent in $G_2$. Recall that $\{u,x,z\}\cap R=\emptyset$ in $G_1$. Then $\{u,x,z\}\cap R^*=\emptyset$ in $G_2$. It is possible that $y\in R^*$. This implies that $d_{G_2}(u)=d_{G_2}(z)=3$. Since $R^*$ is also independent in $G_2$, no two vertices of even degree are adjacent except that the edge $xy$ is allowed. By Lemma 2.1, $G_2$ has a cycle $C_2$ distinct from $C$ such that $C_2$ contains all vertices of $V(C)\setminus R^*$ and the edge $xy$. So $z\in V(C_2)$ as $d_{G_2}(z)=3$.

We can extend the edge set of the cycle $C_2$ containing all vertices of $V(C)\setminus R^*$ and $xy$ in $G_2$ to the edge set of a cycle $C_1$ by possibly adding some elements in $I\cup J\cup D\cup F$. Then $C_1$ contains all vertices of $V(C)\setminus R$ and $xy$ in $G_2$. This implies that $C_1$ in $G_1$ contains all vertices of $C$ and $xy$ except some vertices in $R$ and thus $C_1$ is pleasant in $G_1$. By Claim, $|E(C)|=|E(C_1)|$, $|I_1|=|R_1|$ and $|J_1|=|D_1|=|F_1|=0$. Each vertex in $R_1$ is dominated by one vertex or two vertices in $I_1$ and each vertex in $I_1$ dominates one vertex or two vertices in $R_1$. In particular, the vertex $a$ in $R_1$ is precisely dominated by one vertex in $I_1$ if and only if this vertex in $I_1$ precisely dominates the vertex $a$ in $R_1$. Note that $R_1=\{r_1,\ldots,r_k\}$. Then $I_1=\{u_1,\ldots,u_k\}$ such that $\varphi_1|_{I_1}=R_1$. Let $C=r_1S_1r_2\ldots r_k S_k r_1$, where $S_i=s_{i1}\ldots s_{ij_i}$ for $i\in \{1,\ldots,k\}$. Then $C-R_1=C_1-I_1=S_1\cup \ldots \cup S_k$. Since $G$ is cubic and $xy\in E(C)\cap E(C_1)$, we have $r^+_i\neq r^{-}_{i+1}$ for any $i\in \{1,\ldots,k\}$ and thus $|S_i|\geq 2$.

Now let $G_3$ be a subgraph of $G$ obtained from the graph $C\cup C_1\cup G_1[\{u_1r_1,\ldots,u_kr_k\}]$ by contracting every path $S_i$ into an edge $s_{i1}s_{ij_i}$. Then $G_3$ is a cubic graph and $V(G_3)$ has a partition into sets $A,B$ such that $G[A]$ is a matching $M=\{s_{11}s_{1j_1},\ldots, s_{k1}s_{kj_k}\}$ and $G[B]$ is a matching $M'=\{u_1r_1,\ldots, u_kr_k\}$. Note that both $C$ and $C_1$ contain $xy$ and $z$ with $xy\in E(S_i)$ and $z\in V(S_j)$ for some $i,j\in \{1,\ldots,k\}$. Then $G_3$ has a cycle $C_3=r_1s_{11}s_{1j_1}r_2\ldots r_ks_{k1}s_{kj_k}r_1$ with length $3k$ which contains $M$ and precisely one end of each edge in $M'$. By Theorem 2.2, $G_3$ has a cycle of length larger than $3k$ containing $M=\{s_{11}s_{1j_1},\ldots, s_{k1}s_{kj_k}\}$. Note that $C=r_1S_1r_2\ldots r_kS_kr_1$, where $S_i=s_{i1}\ldots s_{ij_i}$ for $i\in \{1,\ldots,k\}$. Then $G$ has a cycle of length larger than $|\{r_1,\ldots,r_k\}|+|S_1|+\cdots+|S_k|=k+j_1+\ldots+j_k=|E(C)|$ containing $xy$ and $z$, and thus $G$ has a $(x,y)$-path containing $z$ of length larger than $P$, a contradiction. This completes the proof of Theorem \ref{re-main1}.
\q

\section{Chords of longest cycles containing a specified edge in $3$-connected graph with $\delta\geq 4$}

In this section, we prove our second main result, Theorem \ref{main2}. While using the techniques of Thomassen \cite{Th2018}, we also use new ideas to ensure the new longest cycle we find containing the specified edge $e$. For convenience, we restate the theorem here. 

\begin{theo}\label{re-main2}
Let $G$ be a $3$-connected graph with minimum degree at least 4 and $e$ be a specified edge of $G$. Then some longest cycle of $G$ containing $e$ has a chord.
\end{theo}

We will use the following result of Thomassen.

\n{\bf Lemma 4.2} {\rm (Thomassen~\cite{Th2018})} {\em ~If $S$ is an even vertex set in a connected graph $G$, then $G$ has a spanning subgraph $H$ such that every vertex in $S$ has odd degree in $H$ and all other vertices have even degree in $H$.
}

Next, we extend Proposition 1 of Thomassen $\cite{Th2018}$. 

\n{\bf Lemma 4.3} {\em ~Let $C$ be a chordless cycle in a graph $G$ of minimum degree at least 3 and $e\in E(C)$ such that $V(G)-V(C)$ is an independent set in $G$. Then $G$ has a cycle $C_1$ containing $e$ such that either $C_1$ is longer than $C$ or $C_1$ has the same length as $C$ and has a chord. 

Moreover, if $C_1$ containing $e$ has the same length as $C$ and has a chord, then $C_1$ can be chosen such that it has a chord incident with the vertex $a$ of $G_1-V(C)$ satisfying the vertex $a$ has a neighbor vertex on $C$ with degree at least 4, $d_{G_1}(a)=3$ and $N_{G_1}(a)\subseteq V(C_1)$, where $G_1$ is minimal subgraph of $G$ satisfying the conditions of Lemma 4.3.}

\p Let $C$ be a chordless cycle containing $e$ in a graph $G$ of minimum degree at least $3$. We will prove either $G$ has a cycle $C_1$ containing $e$ longer than $C$ or $G$ contains a cycle $C_1$ containing $e$ and has a chord but has the same length as $C$.

Let $G_1$ be a subgraph of $G$ such that $G_1$ is minimal satisfying the conditions of Lemma 4.3 as in $\cite{Th2018}$, i.e., if we delete any edge in $G_1-E(C)$ or any vertex in $G_1-V(C)$, then there is a vertex of degree $2$ in $G_1$.

{\bf Claim 1}{\rm (\cite{Th2018})} {\em If $d_{G_1}(u)=3$ for $u\in V(G_1)-V(C)$, then there is a neighbor of $u$ on $C$ has degree $3$ and if $d_{G_1}(u)\geq 4$ for $u\in V(G_1)-V(C)$, then each neighbor of $u$ on $C$ has degree $3$.}

Let $\{Q_1,\ldots,Q_p\}$ be the set of all components of $G_1-E(C)$. Consider each component $Q_i$ of $G_1-E(C)$. Note that $V(G)-V(C)$ is an independent set in $G$ and $C$ has no chord. Then $Q_i$ is a bipartite graph for $1\leq i\leq p$. Suppose that $V(Q_i)$ is partitioned into two subsets $A_i$ and $B_i$ such that $A_i=V(Q_i)-V(C)$ and $B_i=V(Q_i)\cap V(C)$. Since $\delta(G_1)\geq 3$, we have $|B_i|\geq 3$ and $V(C)=\bigcup\limits_{i=1}^{p} B_i$.

As in $\cite{Th2018}$, we select three vertices $r_i,s_i,t_i$ in $B_i$ such that the number of these vertices with degree $3$ is as large as possible. 

{\bf Claim 2}{\rm (\cite{Th2018})} {\em All of $r_i,s_i,t_i$ have degree 3 in $G_1$ unless $Q_i$ is isomorphic to a bipartite graph $H=G[A_i,B_i]$ with $A_i=\{u_i,v_i\}$ and $B_i=\{r_i,s_i,t_i,w_i\}$, where $d_{G_1}(u_i)=d_{G_1}(v_i)=3$ and $d_{G_1}(w_i)=4$. 
}

Now as in $\cite{Th2018}$, we construct a graph $G^*_1$ from $G_1$ as follows:

(1) We add edges $r_is_i,r_it_i,s_it_i$ for $i\in \{1,\ldots,p\}$ and delete all vertices in $G_1-V(C)$.

(2) We contract all edges each of which has one end vertex of degree two on $C$. 

Then $G^*_1$ is a 4-regular graph on $3p$ vertices and $G^*_1$ has a decomposition into a Hamiltonian cycle and $p$ pairwise vertex disjoint triangles. By Theorem 2.3, we have $\chi(G^*_1)=3$. Since each vertex in a triangle must have different colors, we can assign a vertex of each triangle into color red such that the set of all these red vertices is independent in $G^*_1$. Suppose that $e^*=v_1v_k$ is an edge in $G^*_1$ corresponding to a path $P_1=v_1v_2\ldots v_k$ on $C$ containing $e$. In the following, we assign one end vertex of $e^*$ into color red in $G_1^*$ in order to make sure one end vertex of $e$ is not adjacent to a vertex colored red on $C$ in $G_1$. If $|V(P_1)|\leq 3$ and $e=v_1v_2$ or $|V(P_1)|\geq 4$, then we assign $v_1$ into color red. If $|V(P_1)|=3$ and $e=v_2v_3$, then $e^*=v_1v_3$ and we assign $v_3$ into color red. Since the set of all these red vertices is independent in $G^*_1$, the set of all these red vertices is also independent in $G_1$. We denote $R$ by the set of all red vertices in $G_1$. Then each component $Q_i$ corresponds to a unique red vertex in $\{r_i,s_i,t_i\}$. We rename vertices in $\{r_i,s_i,t_i\}$ such that each red vertex in $Q_i$ is $r_i$.

{\bf Claim 3} {\em One end vertex of $e$ is adjacent to a vertex not in $R$ on $C$ in $G_1$.}

\p  If $|V(P_1)|\leq 3$ and $e=v_1v_2$ or $|V(P_1)|\geq 4$ and $e=v_1v_2$, then $v_1$ is red. Thus another neighbor of $v_1$ distinct from $v_2$ on $C_1^*$ is not red. So this vertex is not in $R$ on $C$ in $G_1$. If $|V(P_1)|=3$ and $e=v_2v_3$, then $e^*=v_1v_3$ and $v_3$ is red. Thus another neighbor of $v_3$ distinct from $v_2$ on $C_1^*$ is not red. So this vertex is not in $R$ on $C$ in $G_1$. If $|V(P_1)|=4$ and $e=v_2v_3$, then $e^*=v_1v_4$ and $v_1$ is red. Thus a neighbor $v_4$ of $v_3$ distinct from $v_2$ on $C_1^*$ is not red. So this vertex is not in $R$ on $C$ in $G_1$. If $|V(P_1)|\geq 5$ and $e=v_2v_3$, then a neighbor $v_4$ of $v_3$ distinct from $v_2$ is not in $C^*_1$. So this vertex is not in $R$ on $C$ in $G_1$. If $|V(P_1)|\geq 4$ and $e=v_iv_{i+1}$ for $i\in \{3,\ldots,k-1\}$, then a neighbor $v_{i-1}$ of $v_i$ distinct from $v_{i+1}$ is not in $C^*_1$. So this vertex is not in $R$ on $C$ in $G_1$. 
\q

Next, we perform a detailed analysis handling $e$, while following the proof process of Proposition 1 in $\cite{Th2018}$.

{\bf Claim 4} {\em If $|A_i|\geq 2$, then every vertex in $A_i$ has degree 3 and has a neighbor vertex on $C$ with degree at least 4 and there is a spanning subgraph $Q'_i$ such that every vertex in $A_i$ has degree 2 in $Q'_i$ and every vertex in $B_i\setminus \{r_i\}$ has odd degree in $Q'_i$. Furthermore, If $|B_i|$ is even, then $r_i$ has degree 1 in $Q'_i$ and there is a vertex $u_i\in A_i$ such that $u_ir_i\in E(Q_i)\cap E(Q'_i)$. If $|B_i|$ is odd, then $r_i$ has degree 0 in $Q'_i$ and there is a vertex $u_i\in A_i$ such that $u_ir_i\in E(Q_i)$ but $u_ir_i\notin E(Q'_i)$.
}

\p Suppose that $|A_i|\geq 2$. Since $A_i=V(Q_i)-V(C)$ and $Q_i$ is a component of $G_1-E(C)$, every vertex in $A_i$ has a neighbor vertex in $C$ with degree at least 4. By Claim 1, every vertex in $A_i$ has degree 3.

If $|B_i|$ is even and $Q_i$ is not isomorphic to $H$, then by Lemma 4.2, $Q_i$ has a spanning subgraph $Q'_i$ such that every vertex in $B_i$ has odd degree in $Q'_i$ and every vertex in $A_i$ have even degree in $Q'_i$. Since $Q_i$ is a bipartite graph and every vertex in $A_i$ has degree 3, every vertex in $A_i$ has degree 2 in $Q'_i$. By Claim 2, all of $r_i,s_i,t_i$ have degree 3 in $G_1$. Note that every vertex in $B_i$ has odd degree in $Q'_i$. So $r_i$ has degree 1 in $Q'_i$ and thus there is a vertex $u_i\in A_i$ such that $u_ir_i\in E(Q_i)\cap E(Q'_i)$.

If $Q_i$ is isomorphic to $H$, then by Lemma 4.2, there is a spanning subgraph $Q'_i$ such that every vertex in $B_i=\{r_i,s_i,t_i,w_i\}$ have degree 1 in $Q'_i$. So we rename vertices in $A_i=\{u_i,v_i\}$ such that $u_ir_i\in E(Q_i)\cap E(Q'_i)$.

If $|B_i|$ is odd, then by Lemma 4.2, $Q_i$ has a spanning subgraph $Q'_i$ such that every vertex in $B_i\setminus \{r_i\}$ has odd degree in $Q'_i$ and every vertex in $A_i\cup \{r_i\}$ have even degree in $Q'_i$. Since $Q_i$ is a bipartite graph and every vertex in $A_i$ has degree 3, every vertex in $A_i$ has degree 2 in $Q'_i$. By Claim 2, all of $r_i,s_i,t_i$ have degree 3 in $G_1$. Note that $r_i$ has even degree in $Q'_i$. So $r_i$ has degree 0 in $Q'_i$ and thus there is a vertex $u_i\in A_i$ such that $u_ir_i\in E(Q_i)$ but $u_ir_i\notin E(Q'_i)$.
\q

If $|A_i|=1$, then we assume that $A_i=\{u_i\}$. By Claim 4, if $|A_i|\geq 2$ and $|B_i|$ is even, there is a vertex $u_i\in A_i$ such that $u_ir_i\in E(Q_i)\cap E(Q'_i)$ and if $|A_i|\geq 2$ and $|B_i|$ is odd, then there is a vertex $u_i\in A_i$ such that $u_ir_i\in E(Q_i)$ but $u_ir_i\notin E(Q'_i)$. So every $u_i\in A_i$ corresponds to each $r_i\in B_i$ in $Q'_i$. Let $I=\{u_1,\ldots,u_p\}$ and $R=\{r_1,\ldots, r_p\}$. Then there is a bijection $\varphi_1: I\rightarrow R$ denoted by $\varphi_1(u_i)=r_i$ for $i\in \{1,2,\ldots,p\}$. Let $J=\bigcup\limits_{i=1}^{p} (A_i-\{u_i\})$ and $Q=\bigcup\limits_{i=1}^{p} (Q_i-\{u_i,r_i\})$. We call each vertex in $I$ as a \emph{pleasant vertex} and each vertex in $J$ as an \emph{unpleasant vertex}.

Now we construct a graph $G_2$ from $G_1$ as follows: 

($1$) If $|A_i|=1$, then we contract the edge $u_ir_i$ to a vertex $r^*_i$

($2$) If $|A_i|\geq 2$, then we contract the edge $u_ir_i$ to a vertex $r^*_i$ and delete all edges from $E(Q_i)\setminus E(Q'_i)$.

Let $R^*=\{r^*_1,r^*_2,\ldots, r^*_p\}$. By Claim 4, if $|A_i|\geq 2$, then every vertex in $A_i$ has degree 2 in $G_2$ and $d_{G_2}(b)$ is odd for any $b\in V(C)\setminus R^*$. Furthermore, if $|B_i|$ is even, then $d_{G_2}(r^*_i)=3$ and thus $|N_{G_2}(r^*_i)\setminus V(C)|=1$ and if  $|B_i|$ is odd, then $d_{G_2}(r^*_i)=4$ and thus $|N_{G_2}(r^*_i)\setminus V(C)|=2$. Note that $\varphi_1: I\rightarrow R$ is a bijection. Then there is a bijection $\varphi_2: R\rightarrow R^*$ denoted by $\varphi_2(r_i)=r^*_i$. 

{\bf Claim 5} {\em $N_{G_2}(a)\subseteq V(C)\setminus R^*$ for any $a\in J$.
}

\p Recall that $J=\bigcup\limits_{i=1}^{p} (A_i-\{u_i\})$. Since $V(G)-V(C)$ is an independent set in G, there are no edges between $I$ and $J$. By the construction of $Q'_i$ and $G_2$, there are no edges between $R^*$ and $J$. So $N_{G_2}(a)\subseteq V(C)\setminus R^*$ for any $a\in J$. 
\q

Let $N_{G_2}(r^*_i)=\{u_{i1},\ldots,u_{i_{t_i}}\}$ for $i\in \{1,\ldots,p\}$. Recall that $R$ is independent in $G_1$. Then $R^*$ is also independent in $G_2$. Since $d_{G_2}(b)$ is odd for any $b\in V(C)\setminus R^*$, $V(G_2)-V(C)$ is independent in $G_2$ and every vertex in $G_2-V(C)$ has even degree 2 in $G_2$, we have no two vertices of even degree are adjacent in $G_2$. By Claim 3, one end vertex of $e$ is adjacent to a vertex not in $R$ on $C$ in $G_1$. Thus one end vertex of $e$ is adjacent to a vertex not in $R^*$ on $C$ in $G_2$. So one end vertex of $e$ is adjacent to a vertex of odd degree on $C$ in $G_2$. By Lemma 2.1, $G_2$ has a cycle $C_2$ distinct from $C$ such that $C_2$ contains all vertices of $V(C)\setminus R^*$ and $e$.

We can extend the edge set of the cycle $C_2$ containing $e$ in $G_2$ to the edge set of a cycle $C_1$ containing $e$ by adding some elements in $I\cup J$. Then $C_1$ contains all vertices of $C$ and $e$ except some vertices in $R$. Without loss of generality, assume that $R_1=\{r_1,\ldots, r_k\}\subseteq R$ is the set of all vertices in $C$ but not in $C_1$ and all vertices in $R_1$ on $C$ appears in the order $r_1,\ldots, r_k$. Since $\delta(G)\geq 3$, it is possible that $r^+_i=r^{-}_{i+1}=v_i$ for some $i\in \{1,\ldots,k\}$. Let $R_2=\{r_i\in R_1~|~r^+_i=r^{-}_{i+1}=v_i\}$. Without loss of generality, assume that $R_2=\{r_1,\ldots, r_q\}\subseteq R_1$. We construct some graphs $G'_1$, $C'$ and $C'_1$ from $G_1$, $C$ and $C_1$ by replacing $v_{i}$ by an edge $v_{i1}v_{i2}$ for $i\in \{1,\ldots, q\}$, respectively. Note that $e\in E(C)\cap E(C_1)$. Then $R_2$ is a proper subset of $R_1$ and $e\neq v_{i1}v_{i2}$ for $i\in \{1,\ldots, q\}$. So $q<k$.

Let $I_1$ and $J_1$ be the set of all pleasant vertices and unpleasant vertices in $C'_1$, respectively. Let $a\in R_1$. Then $C'_1$ must contain two vertices adjacent to $a$ on $C'$. Let $b$ be a neighbor of $a$ on $C'$. Then $C'_1$ contains a path $bdc$, where $d\in I_1\cup J_1$. We say the vertex $d$ dominates $a$ as in \cite{Th2018}. Let $b$ and $b'$ be two neighbors of $a$ on $C'$. Then the vertex $d$ dominates $a$. Possibly, the vertex $d$ dominates some neighbor of $c$ on $C'$. Then any element in $I_1\cup J_1$ dominates at most two vertices in $R_1$. Similarly, $C'_1$ contains a path $b'd'c'$, where $d'\in I_1\cup J_1$. Then the vertex $d'$ dominates $a$. Maybe $a$ is precisely dominated by one element when the vertex $d=d'$ dominates $a$. So each vertex $a$ in $R_1$ is dominated by one element or two elements in $I_1\cup J_1$. Recall that any element in $I_1\cup J_1$ dominates at most two vertices in $R_1$. If follows that $|I_1|+|J_1|\geq |R_1|$.

Note that $E(C'_1)=|E(C')|-2|R_1|+2(|I_1|+|J_1|)$. Since $|I_1|+|J_1|\geq |R_1|$ and $C$ is a longest cycle containing $e$ in $G_1$, we have $C'$ is a longest cycle containing $e$ in $G'_1$ and thus $|E(C')|\geq E(C'_1)=|E(C')|-2|R_1|+2(|I_1|+|J_1|)\geq |E(C')|$. This implies that $|E(C')|=|E(C'_1)|$ and $|I_1|+|J_1|=|R_1|$. Each vertex in $R_1$ is dominated by one vertex or two vertices in $I_1\cup J_1$ and each vertex in $I_1\cup J_1$ dominates one vertex or two vertices in $R_1$. In particular, the vertex $a$ in $R_1$ is precisely dominated by one vertex in $I_1\cup J_1$ if and only if this vertex in $I_1\cup J_1$ precisely dominates the vertex $a$ in $R_1$. Thus $C'-R_1=C'_1-I_1\cup J_1$.

{\bf Case 1} {\em  $J_1=\emptyset$.}

Then $C'-R_1=C'_1-I_1$. Suppose that $C'=r_1S_1r_2\ldots r_k S_kr_1$, where $S_i=s_{i1}\ldots s_{ij_i}$ for $i\in \{1,\ldots,k\}$. Then $C'-R_1=C'_1-I_1=S_1\cup \ldots \cup S_k$. By the construction of $C'$ and $C'_1$, we have $r^+_i\neq r^{-}_{i+1}$ in $C'$ and $C'_1$ for any $i\in \{1,\ldots,k\}$ and thus $|S_i|\geq 2$.

Now let $G_3$ be a graph obtained from the graph $C'\cup C'_1\cup \{u_1r_1,\ldots,u_kr_k\}$ by contracting every path $S_i$ into an edge $s_{i1}s_{ij_i}$. Then $G_3$ is a cubic graph and $V(G_3)$ has a partition into sets $A,B$ such that $G[A]$ is a matching $M=\{s_{11}s_{1j_1},\ldots, s_{k1}s_{kj_k}\}$ and $G[B]$ is a matching $M'=\{u_1r_1,\ldots, u_kr_k\}$.
Recall that $e\in E(C)\cap E(C_1)$, $e\neq v_{i1}v_{i2}$ for $i\in \{1,\ldots, q\}$ and $q<k$. Then $e\in E(S_t)$ for some $t\in \{q+1,\ldots,k\}$. Then $G_3$ has a cycle $C_3=r_1s_{11}s_{1j_1}r_2\ldots r_ks_{k1}s_{kj_k}r_1$ with length $3k$ which contains $M$ and precisely one end of each edge in $M'$. By Theorem 2.2, $G_3$ has a cycle of length larger than $3k$ containing $M=\{s_{11}s_{1j_1},\ldots, s_{k1}s_{kj_k}\}$. Note that $C=r_1S_1r_2\ldots r_kS_kr_1$, where $S_i=s_{i1}\ldots s_{ij_i}$ for $i\in \{1,\ldots,k\}$. Then $G$ has a cycle of length larger than $|\{r_1,\ldots,r_k\}|+|S_1|+\cdots+|S_k|=k+j_1+\ldots+j_k=|E(C')|\geq |E(C)|$ containing $e$.

{\bf Case 2} {\em  $J_1\neq \emptyset$.}

Since $|E(C')|=|E(C'_1)|$, we have $|E(C)|=|E(C_1)|$. Note that $C_1$ is a cycle of $G_1$ containing $e$. Recall that $J_1\subseteq J=\bigcup\limits_{i=1}^{p} (A_i-\{u_i\})$. So $J\neq \emptyset$ and thus $|A_i|\geq 2$ for some $i\in \{1,2,\ldots,p\}$. By Claim 4, every vertex in $A_i$ has degree 3 and has a neighbor vertex on $C$ with degree at least 4. By Claim 5, $N_{G_2}(a)\subseteq V(C)\setminus R^*$ for any $a\in J$. So $d_{G_1}(a)=3$ and $N_{G_1}(a)\subseteq V(C_1)$ for $a\in J_1$. So $C_1$ has a chord incident with this vertex $a$ in $G_1$. Since $G_1$ is minimal subgraph of $G$ satisfying the conditions of Lemma 4.3, we have $C_1$ has a chord incident with this vertex $a$ in $G$. Thus $C_1$ containing $e$ has the same length as $C$ and has a chord. This completes the proof of Lemma 4.3.
\q

\n{\bf Lemma 4.4} {\em ~Let $C$ be a chordless cycle in a $3$-connected graph $G$ of minimum degree at least 4 and $e\in E(C)$. Then $G$ has a cycle $C_1$ containing $e$ such that either $C_1$ is longer than $C$, or $C_1$ has the same length as $C$ and has a chord.
}

\p Let $C$ be a chordless cycle in a $3$-connected graph $G$ of minimum degree at least 4 and $e\in E(C)$. Since $G$ is $3$-connected, the edges between $Q$ and $V(C)$ contain a matching with at least 3 edges for each component $Q$ of $G-V(C)$ with $|V(Q)|\geq 3$. Now we consider each component $Q$ of $G-V(C)$. Then $G'$ can be obtained from $G$ by deleting the edge in $Q$ with $|V(Q)|=2$ and deleting edges between $Q$ and $V(C)$ with $|V(Q)|\geq 3$ such that 

$(i)$ each vertex of $C$ has degree at least 3 in $G'$;

$(ii)$ the edges in $G'$ between $Q$ and $V(C)$ contain a matching with at least 3 edges.

Consider the following 4 properties of a component $Q$ of $G'-V(C)$.

$(1)$ $Q$ has only one vertex and there are only 3 edges between $C$ and $Q$.

$(2)$ the edges in $G'$ between $C$ and $Q$ is a matching with only 3 edges.

$(3)$ There are at least 3 vertices in $V(C)\cap Q$ with degree 3 in $G'$ and if $Q$ has no less than 3 vertices, then the edges in $G'$ between $C$ and $Q$ contain a matching with at least 3 edges. 

$(4)$ there are only 4 edges, say $u_1v_0,u_1v_1,u_2v_2,u_3v_3$ between $C$ and $Q$, where $u_1,u_2,u_3\in V(Q)$, $v_0,v_1,v_2,v_3\in V(C)$, $d_{G'}(v_0)=d_{G'}(v_1)=3$ and $d_{G'}(v_2)=d_{G'}(v_3)\geq 4$. 

Now as in \cite{Th2018}, we say that a component $Q$ of $G'-V(C)$ satisfying a least one of $(1)$, $(2)$ and $(3)$ above is a good component and $(4)$ is a bad component, respectively.

Choose $G'$ such that the number of non-good components is minimal and subject to this, the cardinality of edges between $C$ and $G'-V(C)$ is as small as possible. The proof of Theorem 4 in \cite{Th2018} showed the following claim.

{\bf Claim} {\rm (\cite{Th2018})} {\em Let $Q$ be a component of $G'-V(C)$. Then $Q$ is good.}

If we delete edges between the components $Q$ satisfying $(3)$ to $C$ such that all vertices on $C$ still have degree at least 3, then the proof of Theorem 4 in ~\cite{Th2018} showed the following weaker statement $(3')$ is satisfied as follows:

$(3')$ $Q$ has at least 3 neighbors on $C$ and all neighbors of $Q$ on $C$ have degree only 3.

Now as in \cite{Th2018}, we contract each component $Q$ of $G'-V(C)$ into a vertex. We call the resulting graph $H$. Then $C$ is a chordless cycle in a graph $H$ of minimum degree at least 3 and $e\in E(C)$ such that $V(H)-V(C)$ is an independent set in $H$. Now by Lemma 4.3 with $H$ instead of $G$, $H$ has a cycle $C_1$ containing $e$ such that either $C_1$ is longer than $C$ or $C_1$ has the same length as $C$ and has a chord and $C_1$ can be chosen such that it has a chord incident with the vertex $a$ of $H_1-V(C)$ satisfying the vertex $a$ has a neighbor vertex on $C$ with degree at least 4, $d_{H_1}(a)=3$ and $N_{H_1}(a)\subseteq V(C_1)$, where $H_1$ is minimal subgraph of $H$ satisfying the conditions of Lemma 4.3. Then by the construction of $H$, this vertex $a$ can be extend to a component $Q_a$ which satisfies one of $(1)$, $(2)$, $(3')$ above. Note that $N_{H_1}(a)\subseteq V(C_1)$. If $Q_a$ satisfies $(2)$, then $Q_a$ is a matching with only 3 edges. Thus the edges of $C_1$ cannot form a cycle, a contradiction. Note that the vertex $a$ has a neighbor vertex on $C$ with degree at least 4. If $Q_a$ satisfies $(3')$, then all neighbors of $Q_a$ on $C$ have degree only 3, a contradiction. It follows that $Q_a$ satisfies $(1)$. Hence the chord of $C_1$ in $H$ is also a chord of $C_1$ in $G$. This completes the proof of Lemma 4.4.
\q

From Lemma 4.4, Theorem \ref{re-main2} follows immediately.

\section{Chords of longest cycles containing a specified small set in $3$-connected planar graph} 

In this section, we consider the 3-connected planar case. We will prove the following theorem. This theorem, together with Corollary \ref{main1.1}, will generalize Theorem \ref{theo1.2}.

\begin{theo}\label{re-main3} 
\label{re-main3}
Let $G$ be a $3$-connected planar graph with minimum degree at least $4$. Suppose $S$ is a specified set consisting of either three vertices or an edge plus a vertex. Then every longest cycle of $G$ containing $S$ must have a chord.
\end{theo}

\p Let $G$ be a 3-connected planar graph with minimum degree at least $4$. Suppose $S=\{e,x_3\}$ where $e=x_1x_2$, or $S=\{y_1,y_2,y_3\}$, where $x_1, x_2, x_3, y_1, y_2, y_3$ are vertices of $G$. Since $G$ is a $3$-connected planar graph, there exists a cycle of $G$ containing $S$. Suppose that $C$ is a longest cycle containing $S$. If $S=\{e,x_3\}$, then we can choose a vertex $x_4\neq x_1,x_2,x_3$ on $C$. If $S=\{y_1,y_2,y_3\}$ and $N_C(y_1)\neq \{y_2,y_3\}$, then we can choose a neighbor vertex $y_4$ of $y_1$ on $C$ such that $y_4\neq y_2,y_3$. If $S=\{y_1,y_2,y_3\}$ and $N_C(y_1)=\{y_2,y_3\}$, then we can choose a neighbor vertex $y_5$ of $y_2$ on $C$ such that $y_5\neq y_1,y_3$. In either case, $C$ contains a set $S'$ consisting of an edge along with two other vertices and $S\subseteq S'$. By relabeling, suppose that $S'=\{xy,u,v\}$. Then $C$ is a longest cycle containing $S'$ also. 
We will show that  $C$ must have a chord. Suppose not.
Next we will discuss the following two cases based on whether $C$ is separable in $G$.

{\bf Case 1} {\em  $C$ is not separable.}

Since $\delta(G)\geq 4$, there exists an edge $yz\in E(G)$ and a face $F_1$ of $G$ such that $F_1$ is incident with two adjacent edges $xy,yz$. Let $F_2$ be another face of $G$ incident with $xy$. By Tutte's Lemma (i) and (ii), there is a cycle $C^*$ of $G$ containing $xy$ and $yz$ and $|A(B^*)|\leq 3$ for any bridge $B^*$ of $C^*$. Since $G$ is $3$-connected, we have $2\leq |A(B^*)|\leq 3$ and $|A(B^*)|=2$ if and only if $B^*$ is a chord of $C^*$. Since $\delta(G)\geq 4$, there exists an edge $yw\in \partial(F_2)$. Since $xy,yz\in E(C^*)$, we have $yw\notin E(C^*)$ and thus $yw\in E(B^*)$ for some bridge $B^*$ of $C^*$. So $yw\in E(B^*)\cap \partial(F_2)$. By Tutte's Lemma (iii), $|A(B^*)|\leq 2$ and thus $B^*=yw$ is a chord of $C^*$.

{\bf Claim 1} {\em $V(C)\cap I(B^*)=\emptyset$ for any nonchord bridge $B^*$ of $C^*$.}

\p For any nonchord bridge $B^*$ of $C^*$, $A(B^*)$ is a $3$-vertex-cut of $G$ since $2\leq |A(B^*)|\leq 3$ and $|A(B^*)|=2$ if and only if $B^*$ is a chord of $C^*$. Note that $V(B^*)=I(B^*)\cup A(B^*)$. Consider a partition set $V^*_1$ and $I(B^*)$ of $V(G)-A(B^*)$. Since $C$ is not separable, $V(C)\cap V^*_1=\emptyset$ or $V(C)\cap I(B^*)=\emptyset$. It suffices to show that $V(C)\cap V^*_1\neq \emptyset$. Clearly, $V(C^*)\setminus A(B^*)\subseteq V^*_1$. Suppose that $V(C)\cap V^*_1=\emptyset$. Then $V(C)\cap (V(C^*)\setminus A(B^*))=\emptyset$, i.e, $(V(C)\cap V(C^*))\cap \overline{A(B^*)}=\emptyset$. So $V(C)\cap V(C^*)\subseteq A(B^*)$. Note that $xy\in E(C)\cap E(C^*)$ and $|A(B^*)|=3$. Then either $|V(C)\cap V(C^*)|=2$ or $|V(C)\cap V(C^*)|=3$. Note that $x,y,z$ are three consecutive vertices of $C^*$.

{\bf Subcase 1.1} $|V(C)\cap V(C^*)|=2$. Then $E(C)\cap E(C^*)=\{xy\}$ and $\{x,y\}=V(C)\cap V(C^*)\subseteq A(B^*)$. Note that $yw$ is a chord of $C^*$ and $w\notin \{x,y,z\}$. Since $\{xy,yw\}\subseteq \partial(F_2)$, we have that $yw$ and $I(B^*)$ cannot be both in $int(C^*)$ or both in $ext(C^*)$. So either $yw\in int(C^*)$ and $I(B^*)\subseteq ext(C^*)$ or $yw\in ext(C^*)$ and $I(B^*)\subseteq int(C^*)$. In either case, the closed curve $C':=xCywC^*x$ partitions the interior and exterior of $G$ such that either $xy\in int(C')$ and $yz\in ext(C')$ or $xy\in ext(C')$ and $yz\in int(C')$. By the Jordan Curve Theorem, there is no face incident with both $xy$ and $yz$, a contradiction  to the fact that $\{xy,yz\}\subseteq \partial(F_1)$.

{\bf Subcase 1.2} $|V(C)\cap V(C^*)|=3$ and $|E(C)\cap E(C^*)|=1$. Then $E(C)\cap E(C^*)=\{xy\}$. Assume that $|V(C)\cap V(C^*)|=\{x,y,t\}$ and it is possible that $t\in \{z,w\}$. Note that $yw$ is a chord of $C^*$ and $w\notin \{x,y,z\}$. By a similar argument to that in the Subcase 1.1, we deduce that there is no face incident with both $xy$ and $yz$; a contradiction to the fact that $\{xy,yz\}\subseteq \partial(F_1)$.

{\bf Subcase 1.3} $|V(C)\cap V(C^*)|=3$ and $|E(C)\cap E(C^*)|=2$. Then let $r,x,y,z$ be four consecutive vertices of $C^*$ and thus $E(C)\cap E(C^*)=\{rx,xy\}$ or $E(C)\cap E(C^*)=\{xy,yz\}$. Note that $yw$ is a chord of $C^*$ and $w\notin \{x,y,z\}$. If $E(C)\cap E(C^*)=\{rx,xy\}$, then by a similar argument to that in the Subcase 1.1, we deduce that there is no face incident with both $xy$ and $yz$, a contradiction to the fact that $\{xy,yz\}\subseteq \partial(F_1)$. If $E(C)\cap E(C^*)=\{xy,yz\}$, then $wyxCzC^*w$ is a cycle of $G$ containing $xy$ and all vertices of $C$ larger than $C$, a contradiction.
\q

By Claim 1, $V(C)\subseteq V(C^*)$. Combining this with the fact that $C$ is a longest cycle of $G$ containing $xy$ and $u, v$ and $xy\in E(C)\cap E(C^*)$, we have $V(C)=V(C^*)$. So $G[C]=G[C^*]$. Then $ch(C)+|E(C)|=|E(G[C])|=|E(G[C^*])|=ch(C^*)+|E(C^*)|$. So  $ch(C)=ch(C^*)$. Since $yw$ is a chord of $C^*$, we have $ch(C)=ch(C^*)>1$, a contradiction.

{\bf Case 2} {\em  $C$ is separable.}

In this case, choose a $3$-vertex-cut $T$ of $G$ and a partition set $V_1$ and $V_2$ of $V(G)-T$ with $V(C)\cap V_i\neq \emptyset$ for $i=1,2$ such that $xy\in E(G[T\cup V_1])$ and $|V_2|$ is as small as possible. Since $C$ must pass through at least two vertices of $T$, $C\cap G[T\cup V_i]$ are two paths for $i=1,2$. Let $P_i=C\cap G[T\cup V_i]$ for $i=1,2$ and $T=\{a,b,c\}$. Without loss of generality, assume that $a,b$ are two end vertices of $P_i$. Then $C=aP_2bP_1a$ and $xy\in E(P_1)$ (maybe $xy=ac$ or $xy=bc$ or $xy\notin E(G[T])$). Construct a new graph $G^*$ as follows:
\[G^*:=\left\{
\begin{array}{ll}
   G[T\cup V_2]\cup \{ra,rb,rc\},                               &\mbox{if}~ c\notin V(P_1),                    \\
   G[T\cup V_2]\cup \{ra,rb\}~\mbox{with}~r=c,                  &\mbox{if}~ c\in V(P_1).
\end{array}
\right.\]
Let $C^*=P_2\cup \{ra,rb\}$. Since $C$ is a longest cycle containing $S'=\{xy,u,v\}$, we have $C^*$ is a longest cycle of $G^*$ containing $ra, rb$ and if $f\in \{u,v\}$ is in $V_2$, then $f$ is in $C^*$. Let $F_1$ and $F_2$ be two faces with $ra\in \partial(F_1)\cap \partial(F_2)$. By Tutte's Lemma (i) and (ii), there is a cycle $C'$ of $G^*$ containing $ra, rb$ and $|A(B')|\leq 3$ for any bridge $B'$ of $C'$.

{\bf Claim 2} {\em  $|V(C)\cap V_2|\geq 2$.}

\p Suppose $|V(C)\cap V_2|=1$. Let $V(C)\cap V_2=\{s\}$. If $c\in V(P_1)$, then $c\notin V(P_2)$, and thus $P_2=asb$. If $c\notin V(P_1)$, then $P_2=asb$ or $ascb$ or $acsb$. Since $\delta(G)\geq 4$, there is an edge $st\notin E(P_2)$ with $t\in V_2$. Then $t\in I(B)$ for some bridge $B$ of $C$ and $s\in A(B)$. Since $V(C)\cap V_2=\{s\}$, we have $A(B)\subseteq \{a,b,c,s\}$. Note that $I(B)\neq \emptyset$ and $G$ is 3-connected. Then $|A(B)|\geq 3$. So $|A(B)\cap \{a,b,c\}|\geq 2$. Thus we can find a cycle containing $xy$ and all vertices of $C$ larger than $C$, a contradiction. So $|V(C)\cap V_2|\geq 2$. 
\q

{\bf Claim 3} {\em  $|N_{G}(p)\cap V_2|\geq 2$ for any $p\in T$. }

\p Suppose $|N_{G}(p)\cap V_2|\leq 1$ for some $p\in T$. Then $|N_{G}(p)\cap V_2|=1$ as $G$ is $3$-connected. Let $pp^*\in E(G)$ with $p^*\in V_2$. By Claim 2, $|V(C)\cap V_2|\geq 2$. Then $T^*=(T\setminus \{p\})\cup \{p^*\}$ is also a $3$-vertex-cut of $G$ and there is a partition set $V_1^*=V_1\cup \{p\}$ and $V_2^*=V_2\setminus \{p\}$ of $V(G)-T^*$ with $V(C)\cap V^*_i\neq \emptyset$ for $i=1,2$ such that $xy\in E(G[T^*\cup V^*_1])$ as $xy\in E(G[T\cup V_1])$. But $|V_2^*|<|V_2|$, a contradiction  to the minimality of $|V_2|$. Thus $|N(p)\cap V_2|\geq 2$ for any $p\in T$. 
\q

By Claim 3, $|N_{G}(a)\cap V_2|\geq 2$. Since $\delta(G)\geq 4$, choose two edges $aa_1,aa_2$ with $a_1,a_2\in V_2$ such that $aa_1\in \partial(F_1)$ and $aa_2\in \partial(F_2)$. By Tutte's Lemma (iii), $aa_i$ is either an edge of $C'$ or a chord of $C'$ for $i=1,2$. In either case, $\{a,a_1,a_2\}\subseteq V(C')$. Since $ra\in E(C')$, at least one of $aa_1,aa_2$ is a chord of $C'$.

{\bf Claim 4} {\em $V(C^*)\cap I(B')=\emptyset$ for each nonchord bridge $B'$ of $C'$.}

\p Suppose that $V(C^*)\cap I(B')\neq \emptyset$ for some nonchord bridge $B'$ of $C'$. Then $|I(B')|\geq 1$ and $|A(B')|=3$ as $|A(B')|\leq 3$.

If $r\notin A(B')$, then $c\notin I(B')$ as $r\in V(C')$ and $rc\notin E(C')$. If $r=c$, then $r=c\in V(C')$ and thus $r=c\notin I(B')$. In either case, $\{r,a,b,c\}\cap I(B')=\emptyset$. If $V(C')\setminus \{r,a,b,c\}=\emptyset$. Then $\{a_1,a_2\}=\{b,c\}$. By Claim 2, $V(C)\cap V_2\neq \emptyset$. Then $ab$ is a chord of $C$, a contradiction. Hence $V(C')\setminus \{r,a,b,c\}\neq \emptyset$. Combining this with $\{r,a,b,c\}\cap I(B')=\emptyset$,  we have $I(B')\subseteq V_2$. Clearly, $A(B')$ is also a $3$-vertex-cut of $G$. Since $V(C^*)\cap I(B')\neq \emptyset$, there is a partition set $V_1'=V(G)\setminus (A(B')\cup I(B'))$ and $V_2'=I(B')$ with $V(C)\cap V'_i\neq \emptyset$ for $i=1,2$ such that $xy\in E(G[A(B')\cup V_2'])$ as $xy\in E(G[T\cup V_1])$. But $|V_2'|<|V_2|$, a contradiction  to the minimality of $|V_2|$.

So $r\in A(B')$ and $r\neq c$. Since $ra,rb\in E(C')$ and $|I(B')|\geq 1$, we have $c\in I(B')$ and $rc\in E(B')$. Since $d_{G^*}(r)=3$, we have $rc\in \partial(F_1)\cup \partial(F_2)$. By Tutte's Lemma (iii), $|A(B')|\leq 2$, a contradiction as $|A(B')|=3$.
\q

By Claim 4, $V(C^*)\subseteq V(C')$. Combining this with the fact that $C^*$ is a longest cycle of $G^*$ containing $ra, rb$ and if $f\in \{u,v\}$ is in $V_2$, then $f$ is in $C^*$ and $ra,rb\in E(C^*)\cap E(C')$, we have $V(C^*)=V(C')$. So $ch(C^*)=ch(C')$. Recall that $\{a,a_1,a_2\}\subseteq V(C')$ and at least one of $aa_1,aa_2$ is a chord of $C'$. Then $C^*$ has at least one chord as $ch(C^*)=ch(C')>1$, which implies $C$ has a chord in $G$, a contradiction. This completes the proof of Theorem 5.1.
\q

\section{A new conjecture}

Let $C$ be a circuit of a matroid and $e$ be a non-loop element not in $C$. We say that $e$ is a {\it chord} of $C$ if $e$ is in the closure of $C$. A natural question is whether any largest circuit in a 3-connected matroid has a chord. $AG(3,2)$ shows that the answer is no; indeed, no largest circuit in this matroid has a chord. $M^*(K_n)$ shows that even if we restrict to cographic matroids, the answer is still no.

A graph $F$ is called a linear forest if every component of $F$ is a (possibly empty) path. Clearly, $F=\{e\}$ is a linear forest.  We say that a cycle passes through $F$ if the cycle passes through $E(F)\cup V(F)$ in $G$. In \cite{Hu2001}, Hu, Tian, and Wei proved the following result (we state the slightly weaker version here).

\begin{theo}
Let $G$ be a $k$-connected graph ($k\ge 2$) and let $F$ be a linear forest subgraph of $G$ with $l$ edges and $t$ isolated vertices such that $l+t\le k-2$. Then $G$ has a cycle of length at least min$\{|V(G)|, 2\delta(G)-l\}$ passing through $F$. 
\end{theo}

It is natural to ask whether more general results than our main theorems are true. So we propose the following more general conjecture than conjecture \ref{conj1.4}.

\begin{conj}
Let $G$ be a $k$-connected graph ($k\ge 2$) and let $F$ be a linear forest subgraph of $G$ with $l$ edges and $t$ isolated vertices such that $l+t\le k-2$. Then every longest cycle of $G$ passing through $F$ has a chord.
\end{conj}




\end{document}